\documentclass[a4paper,draft]{amsproc}
\usepackage{amssymb}
\usepackage{amscd} %% Package for commutative diagrams
\theoremstyle{plain}
 \newtheorem{thm}{Theorem}[section]
 \newtheorem{prop}{Proposition}[section]
 \newtheorem{pro}{Problem}[section]

  \newtheorem{conjecture}{Conjecture}[section]
\theoremstyle{definition}
 \newtheorem{exm}{Example}[section]

\theoremstyle{remark}
 
 \numberwithin{equation}{section}

\def\n{\nabla}

\def\f#1#2{\frac{#1}{#2}}

\def\wtd{\widetilde}

\def\a{\alpha}

\def\de{\delta}
\def\De{\Delta}

\def\la{\lambda}

\def\om{\omega}

\def\th{\theta}

\def\w{\wedge}

\def\R{\Bbb{R}}

\def\ra{\rightarrow}

\renewcommand{\leq}{\leqslant}
\renewcommand{\geq}{\geqslant}

\title[]{An intrinsic rigidity theorem for closed minimal hypersurfaces in $\mathbb{S}^5$
with constant nonnegative scalar curvature}



\author[]{Bing Tang    }
\address{
School of Mathematical Sciences \\ % \hfill (Received 00 00 2010)\\
Fudan University   \\ %\hfill (Revised  00 00 2010)\\
Shanghai\\
China}
\email{12210180003@fudan.edu.cn}

\author[]{Ling Yang}
\address{
School of Mathematical Sciences \\ % \hfill (Received 00 00 2010)\\
Fudan University   \\ %\hfill (Revised  00 00 2010)\\
Shanghai\\
China }
\email{yanglingfd@fudan.edu.cn}
\thanks{The authors are partially supported by NSFC} %% optional
\begin{document}

%{\begin{flushleft}\baselineskip9pt\scriptsize
%PUBLICATIONS DE L'INSTITUT MATH\'EMATIQUE\newline
%Nouvelle s\'erie, tome 91(105) (2012), od--do \hfill DOI:
%\end{flushleft}}
\vspace{18mm} \setcounter{page}{1} \thispagestyle{empty}

\begin{abstract}
Let $M^4$ be a closed minimal hypersurface in $\mathbb{S}^5$ with constant nonnegative scalar curvature. Denote by
$f_3$ the sum of the cubes of all principal curvatures, by $g$ the number of distinct principal curvatures.  We prove that, if both $f_3$ and $g$ are constant, then $M^4$ is isoparametric. Moreover, We give all possible values for squared length of the second fundamental form of $M^4$. This result provides another piece of supporting evidence to the Chern conjecture.
\end{abstract}

\maketitle

\section{Introduction}\label{intr}
More than 40 years ago, S.S.Chern proposed the following problem in several places (see \cite{Chern1},\cite{Chern2}):\par
\begin{pro}   Let $M^n$ be a closed minimal submanifold in $\mathbb{S}^{n+m}$ with the second fundamental form of constant length, denote by $\mathcal{A}_n$ the set of all the possible values for the squared length of the second fundamental form of $M^n$, is $\mathcal{A}_n$ a discrete set?
\end{pro}
The affirmative hand of this question is usually called the \textit{Chern conjecture}.

Denote by $B$ the second fundamental form of $M^n$ and let $S:=|B|^2$. Using the Gauss equations, one can easily deduces that
$$S=n(n-1)-R$$
with $R$ denoting the scalar curvature of $M^n$. It means $S$ is in fact an intrinsic geometric quantity, and the Chern conjecture
is equivalent to claiming that the scalar curvature $R$ has gap phenomena for closed minimal submanifolds in Euclidean spheres.

Up to now, it is far from a complete solution of this problem, even in the case that $M$ is a hypersurface (see Problem 105 in \cite{Yau}). Moreover, because all
known examples of closed minimal hypersurfaces in $\Bbb{S}^{n+1}$ with constant scalar curvature are all isoparametric hypersurfaces (the definition
of isoparametric hypersurfaces shall be introduced in Section \ref{iso}). Mathematicians turned
the hypersurface case of Chern conjecture into the following new formulation (see \cite{Ver} \cite{Scherfner2}):

\begin{conjecture}
Let $M^n$ be a closed minimal hypersurface in $\Bbb{S}^{n+1}$ with constant scalar curvature. Then $M$ is an isoparametric hypersurface.
\end{conjecture}

When $n=2$, this conjecture is trivial.
For the case that $n=3$, S. Chang\cite{Chang1,Chang2} gave a positive answer to the Chern conjecture. More precisely, it was shown that any closed minimal hypersurface $M^3$
in $\Bbb{S}^4$ with constant scalar curvature has to be isoparametric, and $\mathcal{A}_3=\{0,3,6\}$.\par
For $n\geq 4$, the Chern conjecture remains open, although some partial result exist for low dimensions and with additional conditions for the curvature functions, such as:

%T. Lusala, M. Scherfner and LAM Souse Jr. gave a partial answer to the 4-dimensional hypersurface case as follows:
\begin{thm}\cite{Lusala1}\label{Thm1}
Let $M^4$ be a closed minimal Willmore hypersurface in $\mathbb{S}^5$ with constant nonnegative scalar curvature. Then $M^4$ is isoparametric.
\end{thm}
%Later, M.Scherfner, L. Vrancken and S.Weiss proved a similar result for hypersurface in $\mathbb{S}^7$:
\begin{thm}\label{Thm2}\cite{Scherfner1}
Let $M^6\in\mathbb{S}^7$ be a closed hypersurface with $H=f_3=f_5\equiv 0$, constant $f_4$ and
$R\geq 0$. Then $M^6$ is isoparametric.
\end{thm}
 Here and in the sequel
$$f_k:=\sum_{i=1}^n \la_i^k$$
with $\la_1,\cdots,\la_n$ being the principal curvatures of $M$.

Note that in Theorem \ref{Thm1}, the 'Willmore' condition equals to saying that $f_3\equiv 0$. It is natural to ask
whether this conclusion holds when '$f_3\equiv 0$' is replaced by a weaker condition that $f_3\equiv const$.
In this paper, we give a partial positive answer to the above question and obtain the main theorem as follows:
\begin{thm}\label{Thm3}
Let $ M^{4}$ be a closed minimal hypersurface in $\mathbb{S}^{5}$ with constant nonnegative scalar curvature. If $f_{3}$ and the number $g$ of distinct principal curvatures of $ M^{4}$ are constant, then $ M^{4}$ is isoparametric.
\end{thm}
Finally, in conjunction with the theory of isoparametric hypersurfaces in Euclidean spheres, we arrive at a classification result (see Theorem \ref{Thm4}), which gave
a piece of supporting evidence to the Chern conjecture.

\bigskip\bigskip

\section{Isoparametric minimal hypersurfaces in $\Bbb{S}^5$}\label{iso}

Let $M^n$ be an immersed hypersurface in $\Bbb{S}^{n+1}$. If $M^n$ has constant principal curvatures, then $M^n$ is said to be an \textit{isoparametric hypersurface}.
Each isoparametric hypersurface is an open subset of a level set of a so-called \textit{isoparametric function} $f$. More precisely, there exists a smooth function $f:\Bbb{S}^{n+1}\ra \R$ and $c\in \R$, such that $|\bar{\n} f|^2$ and $\bar{\De} f$ are both smooth functions of $f$ ($\bar{\n}$ and $\bar{\De}$ are respectively the gradient operator
and Laplace-Beltrami operator on $\Bbb{S}^{n+1}$), and $f(p)=c$ for each $p\in M$. Conversely, given an isoparametric function $f$, the level sets of $f$ consist of
a smooth family of isoparametric hypersurfaces and 2 minimal submanifolds of higher codimension (called \textit{focal submanifolds}).

The following theorem reveals some important geometric properties of isoparametric minimal hypersurfaces in Euclidean spheres (cf. \cite{Cartan1}\cite{Cartan2}\cite{Munzner1}\cite{Munzner2}).

\begin{thm}\label{iso-min}
Let $f:\Bbb{S}^{n+1}\ra \R$ be an isoparametric function, then there exists a unique $c_0\in \R$, such that
$M:=\{x\in \Bbb{S}^{n+1}:f(x)=c_0\}$ is an isoparametric minimal hypersurface. Let $g$ be the number of distinct principal curvatures of $M$, $\la_1>\cdots>\la_g$ be the distinct
principal curvatures, whose multiplicities are $m_1,\cdots,m_g$, respectively, and the denotation of $S$ and $R$ is same as above. Then

\begin{enumerate}
\item $g=1,2,3,4$ or $6$.
\item If $g=1$, then $M$ has to be the totally geodesic great subsphere.
\item If $g=2$, then $M$ has to be a Clifford hypersurface, i.e.
$$M=M_{r,s}:=\Bbb{S}^r\left(\sqrt{\f{r}{n}}\right)\times \Bbb{S}^s\left(\sqrt{\f{s}{n}}\right),$$
where $1\leq r<s\leq n$ and $r+s=n$.
\item If $g=3$, then $m_1=m_2=m_3=2^r$ ($r=0,1,2$ or $3$).
\item There exists $\th_0\in (0,\f{\pi}{g})$, such that
$$\aligned
&\la_k=\cot\left(\f{(k-1)\pi}{g}+\th_0\right),\qquad k=1,\cdots,g,\\
&m_k=m_{k+2} \quad (k\ mod\ g).
\endaligned$$
\item $R\geq 0$ and $S=(g-1)n$.
\end{enumerate}

\end{thm}

E. Cartan \cite{Cartan3} constructed an example of minimal hypersurface in $\mathbb{S}^5$:
\begin{exm}
Denote
$$F:=\left(\sum_i^3(x_i^2-x^2_{i+3})\right)^2+4\left(\sum_i^3x_ix_{i+3}\right)^2$$
For a number $t$ with $0<t<\pi/4$, we denote by $M^4(t)$ a hypersurface in $S^5$ defined by the equation
$$F(x)=\cos^2{2t} ,\qquad x=(x_1,\ldots,x_6)\in\mathbb{S}^5 .$$
A straightforward calculation shows $f:=F|_{\Bbb{S}^5}$ is an isoparametric function and $M^4(\frac{\pi}{8})$ is a minimal isoparametric hypersurface with 4 distinct
principal curvatures,
which is usually called the \textit{Cartan minimal hypersurface}.
\end{exm}

R. Takagi\cite{TAKAGI} proved that $M^4(\frac{\pi}{8})$,
up to congruence, is the unique isoparmetric hypersurface in $\mathbb{S}^5$ with $4$ distinct principal curvatures.
In conjunction with Theorem \ref{iso-min}, we obtain the following result:

\begin{prop}\label{pro}
Let $M^4$ be an isoparametric minimal hypersurface in $\Bbb{S}^5$, then $M^4$, up to a congruence, is either an equator $S^3$, a Clifford hypersurface $(\mathbb{S}^1\left(\frac{1}{2}\right)\times \mathbb{S}^3\left(\frac{\sqrt{3}}{2}\right)$ or
 $\mathbb{S}^2\left(\frac{\sqrt{2}}{2}\right)\times \mathbb{S}^3\left(\frac{\sqrt{2}}{2}\right))$ or then Cartan minimal hypersurface $M^4(\frac{\pi}{8})$,
 and $S=0,4$ or $12$.
\end{prop}

\bigskip\bigskip

\section{Proof of the main theorem }

Let $M^4$ be an immersed hypersurface in $\Bbb{S}^5$. If $\nu$ is a local unit normal vector field along $M$, then 
there exists a pointwise symmetric bilinear form $h$ on $T_p M$, such that
$$B=h\nu.$$
If $\{\omega_1,\omega_2,\omega_3,\omega_4\}$ is a smooth orthonormal coframe field, then $h$ can be written as
$$h=h_{ij}\om_i\otimes \om_j.$$
The covariant derivative $\n h$ with components $h_{ijk}$ is given by
\begin{equation}\label{dh}
\sum_k h_{ijk}\om_k=dh_{ij}+\sum_k h_{kj}\om_{ik}+\sum_k h_{ik}\om_{jk}.
\end{equation}
Here $\{\omega_{ij}\}$ is the connection forms of $M^4$ with respect to $\{\omega_1,\omega_2,\omega_3,\omega_4\}$, which satisfy the following structure equations:
\begin{equation}\label{str}
\aligned
d\om_i&=-\sum_j \om_{ij}\w \om_j,\quad \om_{ij}+\om_{ji}=0\\
d\om_{ij}&=-\sum_k \om_{ik}\w \om_{kj}+\f{1}{2}\sum_{k,l}R_{ijkl}\om_k \w \om_l
\endaligned
\end{equation}
with $R_{ijkl}$ denoting the coefficients of the Riemannian curvature tensor on $M^4$.

In this section, we shall give a proof of the main theorem in Section \ref{intr}.
 
\renewcommand{\proofname}{\bf Proof of Theorem \ref{Thm3}}
\begin{proof}
We shall consider this problem case by case, according to the value of $g$, i.e. the number of distinct principal curvatures.

\textbf{Case I:} $g=1.$

In this case, all the principal curvature are equal to $0$ and hence
 $M^4$ is totally geodesic.

\textbf{Case II:} $g=2$.

 Let $\lambda$ and $\mu$ be distinct pricipal curvatures of $M^4$ with multiplicities $m_1=k, m_2=4-k$, respectively. We need to show that $\lambda, \mu$ are indeed constant functions.\par
Since $\lambda\neq\mu$, from
\begin{equation}\label{g=2}
\begin{aligned}
    m_1\lambda+m_2\mu&=0  \\
    m_1\lambda^2+m_2\mu^2&=S
 \end{aligned}
\end{equation}
we can solve $m_1$, $m_2$ in terms of $\lambda$, $\mu$ and $S$, in other words, $m_1,m_2$ can be seen as continuous functions
of $\lambda,\mu$ and $S$. In conjunction with the fact that $m_1,m_2$ takes values in $\Bbb{Z}$, both $m_1$, $m_2$ are constant, so does $k$. Again from \eqref{g=2}, we have
\begin{equation}\label{}
\displaystyle\lambda=\frac{\sqrt{k(4-k)S}}{2 k}~,~~\displaystyle\mu=-\frac{\sqrt{kS}}{2\sqrt{4-k}},
\end{equation}
or
\begin{equation}\label{}
\displaystyle\lambda=-\frac{\sqrt{k(4-k)S}}{2k}~,~~\displaystyle\mu=\frac{\sqrt{kS}}{2\sqrt{4-k}}.
\end{equation}
Thus $\lambda$ and $\mu$ are both constant and $M^4$ is an isoparametric hypersurface.

\textbf{Case III:} $g=3$.

Let $\lambda,\mu,\sigma$ be distinct principal curvatures of $M^4$, with multiplicities $p,q,r$, respectively, then
\begin{equation}\label{g=3}
\left\{
    \begin{aligned}
    p+q+r&=4\\
    p\lambda+q\mu+r\sigma&=0\\
    p\lambda^2+q\mu^2+r\sigma^2&=S\\
    p\lambda^3+q\mu^3+r\sigma^3&=f_3
    \end{aligned}
\right.
\end{equation}
As in Case II, one can show $p,q,r$ are all constant integer-valued functions. Differentiating both sides of (\ref{g=3}) gives
\begin{equation}\label{}
\left\{
    \begin{aligned}
    pd\lambda+qd\mu+rd\sigma&=0\\
    p\lambda d\lambda+q\mu d\mu+r\sigma d\sigma&=0\\
    p\lambda^2d\lambda+q\mu^2d\mu+r\sigma^2d\sigma&=\frac{1}{3}df_3=0
    \end{aligned}
\right.
\end{equation}
It follows that
\begin{equation}\label{eq}
    \frac{p d\lambda}{\sigma-\mu}=\frac{qd\mu}{\lambda-\sigma}=\frac{rd\sigma}{\mu-\lambda}=\frac{df_3}{3D}=0
\end{equation}
where $D:=(\sigma-\mu)(\sigma-\lambda)(\mu-\lambda)$. Hence $\lambda$, $\mu$ and $\sigma$ are all constant and $M^4$ is isoparametric. (In fact, Theorem \ref{iso-min} shows
there exists no isoparametric minimal hypersurface in $\Bbb{S}^5$ with $g=3$,
so this case cannot occur.)

\textbf{Case IV:} $g=4$.

 let $\lambda_{1}<\lambda_{2}<\lambda_{3}<\lambda_{4}$ be distinct principal curvatures of $M^4$. We say that a coframe field $(U,\omega)$ is \textit{admissible} (see \cite{Scherfner1}) if
\begin{enumerate}
  \item $U$ is an open subset of $M^4$,
  \item $\omega:=\{\omega_1,\omega_2,\omega_3,\omega_4\}$ is a smooth orthonormal coframe field on $U$,
  \item $\omega_1\wedge\omega_2\wedge\omega_3\wedge\omega_4$ is the volume form of $M^4$,
  \item $h=\sum_i\lambda_i\omega_i\otimes\omega_i$.
\end{enumerate}\par
Denote by $F:=\{e_1,e_2,e_3,e_4\}$ the dual frame field of $\omega$. Then it is easily-seen that, $(U,\omega)$ is admissible
if and only if $e_i$ is an unit principal vector associated to $\lambda_i$ for each $1\leq i\leq 4$,
and $\{e_1,e_2,e_3,e_4\}$ is an oriented basis associated to the orientation of $M^4$. Therefore, for
every $p\in M$, there exists an admissible coframe field $(U,\omega)$, such that $p\in U$.

Now we introduce a 3-form on $M^4$: for every admissible coframe field $(U,\omega)$, set
\begin{equation}\label{psi}
\begin{aligned}
    \psi:=&\sum_{1\leq i<j\leq4}(\star(\omega_i\wedge\omega_j))\wedge\omega_{ij},\\
  \end{aligned}
\end{equation}
where $\star$ is the Hodge star operator.
 If $(U,\omega)$ and $(\widetilde{U},\widetilde{\omega})$ are both admissible coframe fields, with $W:=U\cap\widetilde{U}\neq\emptyset$, then on $W$,
$\wtd{\om}_i=\a_i \om_i$ for each $1\leq i\leq 4$, where $\a_i=1$ or $-1$ and $\prod_{i=1}^4 \a_i=1$. Denote by $\{\wtd{\om}_{ij}\}$ the connection form
with respect to $(\wtd{U},\wtd{\om})$, then $\wtd{\om}_{ij}=\a_i\a_j \om_{ij}$ and hence
  $$(\star(\widetilde{\omega}_i\wedge\widetilde{\omega}_j))\wedge\widetilde{\omega}_{ij}=(\star(\omega_i\wedge\omega_j))\wedge\omega_{ij}$$
 holds for any $i<j$. Therefore $\psi$ is well-defined on $M^4$.\par

Now we compute the exterior differential of the form $\psi$. Due to the definition of the Hodge star operator, $\psi$ can be written as
\begin{equation}\label{psi2}
\aligned
\psi=&\om_1\w \om_2\w \om_{34}+\om_2\w \om_3\w \om_{14}+\om_3\w \om_1\w \om_{24}\\
&+\om_1\w\om_4\w \om_{23}+\om_2\w \om_4\w \om_{31}+\om_3\w \om_4\w
\om_{12}.
\endaligned
\end{equation}

Substituting $h_{ij}=\la_i \de_{ij}$ into (\ref{dh}), we have
\begin{equation}\label{om}
\om_{ij}=\f{1}{\la_j-\la_i}\sum_{k}h_{ijk}\om_k\qquad \forall i\neq
j.
\end{equation}
Combining (\ref{om}) and (\ref{str}) yields
$$\aligned
d\om_1=&-(\om_{12}\w \om_2+\om_{13}\w \om_3+\om_{14}\w \om_4)\\
=&(\cdots)\w \om_2-\f{1}{\la_3-\la_1}(h_{131}\om_1+h_{134}\om_4)\w \om_3\\
&-\f{1}{\la_4-\la_1}(h_{141}\om_1+h_{143}\om_3)\w \om_4.
\endaligned$$
Hence
\begin{equation}\label{d1}
\aligned
&d\om_1\w \om_2\w
\om_{34}\\
=&-\Bigg[\f{h_{113}h_{443}}{(\la_3-\la_1)(\la_3-\la_4)}+\f{h_{114}h_{334}}{(\la_4-\la_1)(\la_4-\la_3)}+\f{h_{134}^2}{(\la_3-\la_1)(\la_4-\la_1)}\Bigg]\star 1
\endaligned
\end{equation}
(where we have used Codazzi equations). A similar calculation shows
\begin{equation}\label{d2}\aligned
&\om_1\w d\om_2\w \om_{34}\\
=&\Bigg[\f{h_{223}h_{443}}{(\la_3-\la_2)(\la_3-\la_4)}+\f{h_{224}h_{334}}{(\la_4-\la_2)(\la_4-\la_3)}+\f{h_{234}^2}{(\la_3-\la_2)(\la_4-\la_2)}\Bigg]\star
1.
\endaligned\end{equation}
By the structure equations,
\begin{equation}\label{d3}\aligned
d\om_{34}=&-\om_{31}\w \om_{32}\w_{24}+\f{1}{2}\sum_{k,l}R_{34kl}\om_k\w \om_l\\
=&\Bigg[\f{h_{331}h_{441}}{(\la_3-\la_1)(\la_4-\la_1)}+\f{h_{332}h_{442}}{(\la_3-\la_2)(\la_4-\la_2)}-\f{h_{134}^2}{(\la_3-\la_1)(\la_4-\la_1)}\\
 &-\f{h_{234}^2}{(\la_3-\la_2)(\la_4-\la_2)}+R_{3434}\Bigg]\om_3\w \om_4+(\cdots)\w \om_1+(\cdots)\w \om_2.
 \endaligned\end{equation}
Combining (\ref{d1})-(\ref{d3}) gives
\begin{equation}
\aligned
&d(\om_1\w \om_2\w \om_{34})=d\om_1\w \om_2\w \om_{34}-\om_1\w d\om_2\w \om_{34}+\om_1\w \om_2\w d\om_{34}\\
=&\Bigg[\f{h_{331}h_{441}}{(\la_3-\la_1)(\la_4-\la_1)}+\f{h_{332}h_{442}}{(\la_3-\la_2)(\la_4-\la_2)}-\f{h_{113}h_{443}}{(\la_3-\la_1)(\la_3-\la_4)}\\
&-\f{h_{114}h_{334}}{(\la_4-\la_1)(\la_4-\la_3)}-\f{h_{223}h_{443}}{(\la_3-\la_2)(\la_3-\la_4)}+\f{h_{224}h_{334}}{(\la_4-\la_2)(\la_4-\la_3)}\\
&-\f{2h_{134}^2}{(\la_3-\la_1)(\la_4-\la_1)}-\f{2h_{234}^2}{(\la_3-\la_2)(\la_4-\la_2)}+R_{3434}\Bigg]\star 1.
\endaligned
\end{equation}
Similarly, one can compute the exterior differential of each term of (\ref{psi2}); taking the sum of these equations, we arrive at
\begin{equation}\label{3form}
    d\psi=\left(\frac{1}{2}R-\sum_{l=1}^{4}I_{l}\right)\star1.
\end{equation}
where
\begin{equation}\label{B}
    I_{l}=\sum_{l\neq i<j\neq l}\frac{h_{iil}h_{jjl}}{(\lambda_{l}-\lambda_{i})(\lambda_{l}-\lambda_{j})},~~~~\forall l=1,2,3,4.
\end{equation}\par
Taking the exterior differential of
\begin{equation}
\left\{
\begin{aligned}
    \sum_{i}h_{ii}=&0\\
    \sum_{i,j}h_{ij}^2=&S=const\\
    \sum_{i,j,k}h_{ij}h_{jk}h_{ki}=&f_{3}=const\\
 \end{aligned}
 \right.
\end{equation}
implies that
\begin{equation}
\left\{
\begin{aligned}
    \sum_{i}h_{iik}&=0\\
    \sum_i\lambda_{i}h_{iik}&=0\\
    \sum_i\lambda_{i}^{2}h_{iik}&=0\\
 \end{aligned}
\right.
\end{equation}
holds for each $1\leq k\leq 4$. Especially, letting $k:=1$ gives
\begin{equation}
\left\{
\begin{aligned}
    h_{111}+h_{221}+h_{331}+h_{441}&=0\\
    \lambda_{1}h_{111}+\lambda_{2}h_{221}+\lambda_{3}h_{331}+\lambda_{4}h_{441}&=0\\
    \lambda_{1}^{2}h_{111}+\lambda_{2}^{2}h_{221}+\lambda_{3}^{2}h_{331}+\lambda_{4}^{2}h_{441}&=0\\
 \end{aligned}
\right.
\end{equation}
Since $\lambda_{1},\lambda_{2},\lambda_{3}$ and $\lambda_{4}$ are distinct at every point, we can express $h_{ii1}, i=2,3,4$, in terms of $h_{111}$:
\begin{equation}
 h_{ii1}=-\frac{\prod\limits_{j\neq i,1}(\lambda_{j}-\lambda_{1})}{\prod\limits_{j\neq i,1}(\lambda_{j}-\lambda_{i})}h_{111},~~~~\forall i=2,3,4.
\end{equation}
Let $K:=\det h$ be the Gauss-Kronecker curvature of $M^4$ and denote
$$dK=\sum_i K_i \om_i,$$
then
\begin{equation}
K_{1}=\sum_{i=1}^{4}\left(h_{ii1}\prod_{j\neq i}\la_j\right)
=-(\lambda_{1}-\lambda_{2})(\lambda_{1}-\lambda_{3})(\lambda_{1}-\lambda_{4})h_{111}
\end{equation}
and hence
\begin{equation}
h_{ii1}=\frac{K_{1}}{\prod\limits_{j\neq i}(\lambda_{j}-\lambda_{i})}.
\end{equation}
In a similar way, we have
\begin{equation}\label{A}
    h_{iil}=\frac{K_{l}}{\prod\limits_{j\neq i}(\lambda_{j}-\lambda_{i})},~~~~\forall i,l=1,2,3,4.
\end{equation}
Substituting \eqref{A} into \eqref{B}, we deduce that
\begin{equation}\label{C}
I_{l}=K_{l}^{2}\sum_{l\neq i<j\neq l}
\frac{1}{(\lambda_{l}-\lambda_{i})(\lambda_{l}-\lambda_{j})
\prod\limits_{m\neq i}(\lambda_{m}-\lambda_{i})\prod\limits_{m\neq j}(\lambda_{m}-\lambda_{j})}.
\end{equation}
More precisely,
\begin{equation}
\begin{aligned}
 I_{1}=&K_{1}^{2}\sum_{1\neq i<j\neq 1}
\frac{1}{(\lambda_{1}-\lambda_{i})(\lambda_{1}-\lambda_{j})
\prod\limits_{m\neq i}(\lambda_{m}-\lambda_{i})\prod\limits_{l\neq j}(\lambda_{l}-\lambda_{j})}\\
=
&K_{1}^{2}\left[
\frac{1}{(\lambda_{1}-\lambda_{2})(\lambda_{1}-\lambda_{3})\prod_{m\neq2}(\lambda_m-\lambda_2)\prod_{l\neq3}(\lambda_l-\lambda_3)}\right.\\
&~~~~~~~~~+\frac{1}{(\lambda_{1}-\lambda_{2})(\lambda_{1}-\lambda_{4})\prod_{m \neq2}(\lambda_m-\lambda_2)\prod_{l \neq3}(\lambda_l-\lambda_4)}\\
&~~~~~~~~~+\left.\frac{1}{(\lambda_{1}-\lambda_{3})(\lambda_{1}-\lambda_{4})\prod_{m \neq3}(\lambda_m-\lambda_3)\prod_{l \neq4}(\lambda_l-\lambda_4)}\right]\\
=&-\frac{K_{1}^{2}}{D^{2}}[(\lambda_{4}-\lambda_{3})(\lambda_{4}-\lambda_{2})(\lambda_{4}-\lambda_{1})^{2}%%%%%%%%%%%%%%
+(\lambda_{3}-\lambda_{4})(\lambda_{3}-\lambda_{2})(\lambda_{3}-\lambda_{1})^{2}\\
&~~~~~~~~~~~~+(\lambda_{2}-\lambda_{4})(\lambda_{2}-\lambda_{3})(\lambda_{2}-\lambda_{1})^{2}]
 \end{aligned}
\end{equation}
where $\displaystyle D:=\prod_{1\leq i<j\leq4}(\lambda_{j}-\lambda_{i})$. Similarly, one computes
\begin{equation}
\begin{aligned}
 I_{2}=&-\frac{K_{2}^{2}}{D^{2}}[(\lambda_{4}-\lambda_{3})(\lambda_{4}-\lambda_{2})^{2}(\lambda_{4}-\lambda_{1})%%%%%%%%%%%%%%
+(\lambda_{3}-\lambda_{4})(\lambda_{3}-\lambda_{2})^{2}(\lambda_{3}-\lambda_{1})\\
&~~~~~~~~~~~~+(\lambda_{1}-\lambda_{4})(\lambda_{1}-\lambda_{3})(\lambda_{1}-\lambda_{2})^{2})] ,
 \end{aligned}
\end{equation}
 \begin{equation}
\begin{aligned}
 I_{3}=&-\frac{K_{3}^{2}}{D^{2}}[(\lambda_{4}-\lambda_{3})^{2}(\lambda_{4}-\lambda_{2})(\lambda_{4}-\lambda_{1})%%%%%%%%%%%%%%
+(\lambda_{2}-\lambda_{4})(\lambda_{2}-\lambda_{3})^{2}(\lambda_{2}-\lambda_{1})\\
&~~~~~~~~~~~~+(\lambda_{1}-\lambda_{4})(\lambda_{1}-\lambda_{3})^{2}(\lambda_{1}-\lambda_{2})]
 \end{aligned}
\end{equation}
and
\begin{equation}
\begin{aligned}
 I_{4}=&-\frac{K_{4}^{2}}{D^{2}}[(\lambda_{3}-\lambda_{4})^{2}(\lambda_{3}-\lambda_{2})(\lambda_{3}-\lambda_{1})%%%%%%%%%%%%%%
+(\lambda_{2}-\lambda_{4})^{2}(\lambda_{2}-\lambda_{3})(\lambda_{2}-\lambda_{1})\\
&~~~~~~~~~~~~+(\lambda_{1}-\lambda_{4})^{2}(\lambda_{1}-\lambda_{3})(\lambda_{1}-\lambda_{2})] .
 \end{aligned}
\end{equation}
Observing that $\lambda_{1}<\lambda_{2}<\lambda_{3}<\lambda_{4}$, we can derive estimates as follows.
\begin{equation}
\begin{aligned}
I_{1}=&-\frac{K_{1}^{2}}{D^{2}}[(\lambda_{4}-\lambda_{3})(\lambda_{4}-\lambda_{2})(\lambda_{4}-\lambda_{1})^{2}%%%%%%%%%%%%%%
+(\lambda_{3}-\lambda_{4})(\lambda_{3}-\lambda_{2})(\lambda_{3}-\lambda_{1})^{2}\\
&~~~~~~~~~~~~+(\lambda_{2}-\lambda_{4})(\lambda_{2}-\lambda_{3})(\lambda_{2}-\lambda_{1})^{2}]\\
\leq&-\frac{K_{1}^{2}}{D^{2}}[(\lambda_{4}-\lambda_{3})(\lambda_{3}-\lambda_{2})(\lambda_{4}-\lambda_{1})^{2}%%%%%%%%%%%%%%
+(\lambda_{3}-\lambda_{4})(\lambda_{3}-\lambda_{2})(\lambda_{3}-\lambda_{1})^{2}\\
&~~~~~~~~~~~~+(\lambda_{2}-\lambda_{4})(\lambda_{2}-\lambda_{3})(\lambda_{2}-\lambda_{1})^{2}]\\
=&-\frac{K_{1}^{2}}{D^{2}}[(\lambda_{4}-\lambda_{3})(\lambda_{3}-\lambda_{2})    (\lambda_{4}-\lambda_{3})(\lambda_{4}+\lambda_{3}-2\lambda_{1})\\
&~~~~~~~~~~~~+(\lambda_{2}-\lambda_{4})(\lambda_{2}-\lambda_{3})(\lambda_{2}-\lambda_{1})^{2}]\\
\leq& 0.
 \end{aligned}
\end{equation}
\begin{equation}
\begin{aligned}
 I_{2}=&-\frac{K_{2}^{2}}{D^{2}}[(\lambda_{4}-\lambda_{3})(\lambda_{4}-\lambda_{2})^{2}(\lambda_{4}-\lambda_{1})%%%%%%%%%%%%%%
+(\lambda_{3}-\lambda_{4})(\lambda_{3}-\lambda_{2})^{2}(\lambda_{3}-\lambda_{1})\\
&~~~~~~~~~~~~+(\lambda_{1}-\lambda_{4})(\lambda_{1}-\lambda_{3})(\lambda_{1}-\lambda_{2})^{2})]\\
\leq&-\frac{K_{2}^{2}}{D^{2}}[(\lambda_{4}-\lambda_{3})(\lambda_{4}-\lambda_{2})^{2}(\lambda_{4}-\lambda_{1})%%%%%%%%%%%%%%
+(\lambda_{3}-\lambda_{4})(\lambda_{3}-\lambda_{2})^{2}(\lambda_{4}-\lambda_{1})\\
&~~~~~~~~~~~~+(\lambda_{1}-\lambda_{4})(\lambda_{1}-\lambda_{3})(\lambda_{1}-\lambda_{2})^{2})]\\
=&-\frac{K_{2}^{2}}{D^{2}}[(\lambda_{4}-\lambda_{3})(\lambda_{4}-\lambda_{1})(\lambda_{4}-\lambda_{3})(\lambda_{4}+\lambda_{3}-2\lambda_{2})\\
&~~~~~~~~~~~~+(\lambda_{1}-\lambda_{4})(\lambda_{1}-\lambda_{3})(\lambda_{1}-\lambda_{2})^{2})]\\
\leq& 0.
 \end{aligned}
\end{equation}
In the same way, $I_{3}\leq 0 ,I_{4}\leq 0$.\par
Note that $M^4$ is closed. Integrating both sides of \eqref{3form} on $M^4$ and then using Stokes's theorem gives
\begin{equation}
\displaystyle0=\int_{M^{4}}d\psi=\frac{1}{2}\int_{M^{4}}R~\star1-\int_{M^{4}}\sum_{k}I_{k}~\star1 .
\end{equation}
Since $R\geq 0$ and $I_{k}\leq0$ for $k=1,2,3,4$, it follows that $R=0$ and $I_{k}=0 , k=1,2,3,4$. From \eqref{C}, $dK=0$, so $\prod_{i=1}^4 \la_i=K=const$.
In conjunction with $\sum_i \la_i=0$, $\sum_i \la_i^2=S=const$ and $\sum_i \la_i^3=f_3=const$, one can easily deduce that $\la_i$ ($1\leq i\leq 4$) are
all constant on $M$.
Thus $M^4$ is an isoparametric hypersurface.
\end{proof}

Combining Theorem \ref{Thm3} and Proposition \ref{pro} yields a classification theorem as follows.

\begin{thm}\label{Thm4}
Let $ M^{4}$ be a closed minimal hypersurface in $\mathbb{S}^{5}$ with constant nonnegative scalar curvature. If $f_{3}$ and the number $g$ of distinct principal curvatures of $ M^{4}$ are constant, then $M^4$, up to a congruence, is either an equator $S^3$, a Clifford hypersurface $(\mathbb{S}^1\left(\frac{1}{2}\right)\times \mathbb{S}^3\left(\frac{\sqrt{3}}{2}\right)$ or
 $\mathbb{S}^2\left(\frac{\sqrt{2}}{2}\right)\times \mathbb{S}^3\left(\frac{\sqrt{2}}{2}\right))$ or then Cartan minimal hypersurface $M^4(\frac{\pi}{8})$. Let $S$ denote the squared length of the second fundamental form of $M^4$, then $S=0,4$ or $12$.
\end{thm}

%%%%%%%%%%%%%%%%%%%%%%%%%%%%%%
\bigskip

\bibliographystyle{amsplain}

\end{document}